\input ./style/arxiv-general.cfg
\documentclass[MSNbibl,number,citesort,seceqn,secthm,dvips]{arxbj}
\makeatletter
   \@ifpackageloaded{graphicx}{}{\usepackage{graphicx}}
\makeatother
\usepackage{mathrsfs}

\volume{22}
\issue{4}
\pubyear{2016}
\firstpage{2301}
\lastpage{2324}
\doi{10.3150/15-BEJ729}% Updated by VTEXPTS2LaTeX.exe, 14.08.2015 10:34
\docsubty{FLA}

\makeatletter
\newcommand{\defas}{:=}
\newcommand{\CombStruct}[1]{[[#1]]}
\newcommand{\process}[2]{\nprocess{#1}{#2}n}
\newcommand{\nprocess}[3]{(#1_{#3})_{#3 \in #2}}
\newcommand{\Atoms}{\mathscr{A}}
\newcommand{\theset}[1]{\lbrace #1 \rbrace}
\newcommand{\betadist}{\mathrm{beta}}
\newcommand{\poisson}{\operatorname{Poisson}}
\newcommand{\nbdist}{\operatorname{NB}}
\newcommand{\bnbdist}{\mathrm{beta}\mbox{-}\mathrm{NB}}
\newcommand{\digammadist}{\operatorname{digamma}}
\newcommand{\gammadist}{\operatorname{gamma}}

\newcommand{\NBPLAW}{\mathrm{NBP}}
\newcommand{\BNBPLAW}{\mathrm{BNBP}}

\newcommand{\harmonicsum}[1]{\lambda_{#1}}
\newcommand{\Histories}{\mathcal{H}}
\newcommand{\historycount}[1]{M_{#1}}
\newcommand{\historysum}[1]{S({#1})}
\newcommand{\countspace}{{\Histories}_\infty}
\newcommand{\gprocess}[2]{(#1)_{#2}}
\newcommand{\extendedspace}{\mathcal{H}^{(0)}}
\newcommand{\existingspace}{\mathcal{H}^+}
\newcommand{\basicbb}{b}

\newcommand{\cardsharp}{\#}
\newcommand{\bspace}{\Omega}
\newcommand{\bsa}{\mathcal A}
\newcommand{\borelspace}{(\bspace,\bsa)}
\newcommand{\Mspace}{\mathcal M \borelspace}

\newcommand{\gprocessM}[1]{\gprocess{\historycount{h}}{h\in\Histories_{#1}}}

\newcommand{\Xseqn}[1]{X_{[#1]}}
\newcommand{\betafn}{\mathcal{B}}
\newcommand{\oh}{\hslash}
\newcommand{\extendedhist}[2]{\Histories_{#1}^{(#2)}}
\newproclaim{definition}{Definition}[section]
\newremark{remark}{Remark}[section]
\makeatother

\begin{document}
\begin{frontmatter}

\title{The combinatorial structure of beta negative binomial processes}
\runtitle{The combinatorial structure of beta negative binomial processes}

\begin{aug}
%%%% inicialai - be tarpu
% Corresponding author: Creighton Heaukulani - ckheau@gmail.com% Updated by VTEXPTS2LaTeX.exe, 14.08.2015 15:28
%Updated by VTEXPTS2LaTeX.exe, 14.08.2015 10:34
\author[A]{\inits{C.}\fnms{Creighton}~\snm{Heaukulani}\corref{}\thanksref{A}\ead[label=e1]{ckh28@cam.ac.uk}}
%%,
\and
\author[B]{\inits{D.M.}\fnms{Daniel M.}~\snm{Roy}\thanksref{B}\ead[label=e2]{droy@utstat.toronto.edu}}
%
%\author[]{\inits{}\fnms{}~\snm{}\thanksref{}\ead[label=]{}}
%%\runauthor{} %% auto
%\dedicated{}
\address[A]{Department of Engineering, University of Cambridge,
Trumpington Street, Cambridge, CB2 1PZ, United Kingdom.
\printead{e1}}
\address[B]{Department of Statistical Sciences, University of Toronto,
100 St. George Street, Toronto, ON, M5S 3G3, Canada.
\printead{e2}}
\end{aug}

% HISTORY:
%
\received{\smonth{6} \syear{2014}}% Updated by VTEXPTS2LaTeX.exe,
%14.08.2015 10:34
%
\revised{\smonth{3} \syear{2015}}% Updated by VTEXPTS2LaTeX.exe,
%14.08.2015 10:34

% ABSTRACT
%
\begin{abstract}
We characterize the combinatorial structure of conditionally-i.i.d.
sequences of negative binomial processes with a common beta process
base measure.
In Bayesian nonparametric applications, such processes have served as
models for latent \emph{multisets} of features underlying data.
Analogously, random \emph{subsets} arise from conditionally-i.i.d.
sequences of Bernoulli processes with a common beta process base
measure, in which case the combinatorial structure is described by the
\emph{Indian buffet process}.
Our results give a count analogue of the Indian buffet process, which
we call a \emph{negative binomial Indian buffet process}.
As an intermediate step toward this goal, we provide a construction for
the beta negative binomial process that avoids a representation of the
underlying beta process base measure.
We describe the key Markov kernels needed to use a NB-IBP
representation in a Markov Chain Monte Carlo algorithm targeting a
posterior distribution.
\end{abstract}

% KEYWORDS
% visi is mazosios raides ir pagal abecele
%
\begin{keyword}
\kwd{Bayesian nonparametrics}
\kwd{Indian buffet process}
\kwd{latent feature models}
\kwd{multisets}
\end{keyword}
\end{frontmatter}

%s1 #&#
\section{Introduction}\label{sec1}

The focus of this article is on exchangeable sequences of \emph
{multisets}, that is, set-like objects in which repetition is allowed.
Let $\Omega$ be a complete, separable metric space equipped with its
Borel $\sigma$-algebra $\mathcal{A}$ and let $\mathbb{Z}_+\defas\{
0, 1, 2,\ldots\}$ denote the non-negative integers.
By a \emph{point process} on $(\Omega,\mathcal{A})$, we mean a
random measure
$X$ on $(\Omega,\mathcal{A})$ such that $X(A)$ is a $\mathbb
{Z}_+$-valued random
variable for every $A \in\mathcal{A}$.
Because $(\Omega,\mathcal{A})$ is Borel, we may write
%
%e1 #&#
\begin{equation}
\label{eqpointprocess} X = \sum_{k\le\kappa}\delta_{\gamma_k}
\end{equation}
for a random element $\kappa$ in $\overline{\mathbb{Z}}_+\defas
\mathbb{Z}_+\cup\{
\infty\}$ and some -- not necessarily distinct -- random elements
$\gamma
_1, \gamma_2, \ldots$ in $\Omega$.
We will take the point process $X$ to represent the \emph{multiset} of
its unique atoms $\gamma_k$ with corresponding multiplicities $X
\theset
{\gamma_k}$. We say $X$ is \emph{simple} when $X\{\gamma_k\} = 1$ for
all $k \le\kappa$, in which case $X$ represents a set.

In statistical applications, \emph{latent feature models} associate
each data point $y_n$ in a dataset with a latent point process $X_n$
from an exchangeable sequence of simple point processes, which we
denote by $\process{X}{\mathbb{N}}\defas(X_1, X_2, \ldots)$.
The unique atoms among the sequence $\process X {\mathbb{N}}$
are
referred to
as \emph{features},
and a data point is said to \emph{possess} those features appearing in
its associated point process.
We can also view these latent feature models as
generalizations of mixture models that allow
data points to belong to multiple, potentially overlapping clusters
\cite{GG06,broderick2013cluster}.
For example, in an object recognition task, a model for a dataset
consisting of street camera images could associate each image with a
subset of object classes -- for example, ``trees'', ``cars'', and ``houses'',
etc.~-- appearing in the images.
In a document modeling task, a model for a dataset of news articles
could associate each document with a subset of topics -- for example,
``politics'', ``Europe'', and ``economics'', etc. -- discussed in the documents.
Recent work in Bayesian nonparametrics utilizing exchangeable sequences
of simple point processes have focused on the Indian buffet process
(IBP) \cite{GG06,GGS2007}, which characterizes the marginal
distribution of the sequence $\process X {\mathbb{N}}$ when they are
conditionally-i.i.d.  Bernoulli processes, given a common beta process
base measure \cite{Hjort1990,TJ2007}.

If the point processes $\process X {\mathbb{N}}$ are no longer
constrained to
be simple, then data points can contain multiple copies of features.
For example, in the object recognition task, an image could be
associated with two cars, two trees, and one house. In the document
modeling task, an article could be associated with 100 words from the
politics topic, 200 words from the Europe topic, and 40 words from the
economics topic.
In this article, we describe a count analogue of the IBP called the
\emph{negative binomial Indian buffet processes} (NB-IBP), which
characterizes the marginal distribution of $\process X {\mathbb{N}}$
when it
is a conditionally i.i.d. sequence of negative binomial processes
\cite{BMPJpre,ZHDC12}, given a common beta process base measure.
This characterization allows us to describe new Markov Chain Monte
Carlo algorithms for posterior inference that do not require numerical
integrations over representations of the underlying beta process.

%s1.1 #&#
\subsection{Results}

Let $c>0$, let $\widetilde{B}_0$ be a non-atomic, finite measure on
$\Omega$, and
let $\Pi$ be a Poisson (point) process on $\Omega\times(0,1]$
with intensity
%
%e2 #&#
\begin{eqnarray}
\label{eqfirstBPintensity} (\mathrm{d}s, \mathrm{d}p) \mapsto c p^{-1} (1 -
p)^{c -1}\, \mathrm{d}p \widetilde {B}_0(\mathrm{d}s).
\end{eqnarray}
As this intensity is non-atomic and merely $\sigma$-finite, $\Pi$
will have an infinite number of atoms almost surely (\mbox{a.s.}), and
so we
may write $\Pi= \sum_{j=1}^\infty\delta_{(\gamma_j, \basicbb
_j)}$\vspace*{1pt} for some \mbox{a.s.} unique random elements $\basicbb_1,
\basicbb_2, \ldots$ in $(0,1]$ and $\gamma_1, \gamma_2, \ldots$ in $\Omega$.
From $\Pi$, construct the random measure
%
%e3 #&#
\begin{eqnarray}
&& B \defas\sum_{j=1}^\infty
\basicbb_j \delta_{\gamma_j},
\end{eqnarray}
which is a \emph{beta process} \cite{Hjort1990}.
The construction of $B$ ensures that the random variables $B(A_1),
\ldots, B(A_k)$ are independent for every finite, disjoint collection
$A_1,\ldots,A_k\in\mathcal{A}$, and $B$ is said to be \emph
{completely random} or equivalently, have \emph{independent
increments} \cite{Kingman1967}. We review completely random measures
in Section~\ref{secprelims}.

The conjugacy of the family of beta distributions with various other
exponential families carries over to beta processes and randomizations
by probability kernels lying in these same exponential families. The
beta process is therefore a convenient choice for further
randomizations, or in the language of Bayesian nonparametrics, as a
prior stochastic process.
For example, previous work has focused on the (simple) point process
that takes each atom $\gamma_j$ with probability $\basicbb_j$ for every
$j\ge1$, which is, conditioned on $B$, called a \emph{Bernoulli
process} (with base measure $B$) \cite{TJ2007}.
In this article, we study the point process\vspace*{-6pt}
%
%e4 #&#
\begin{eqnarray}
\label{eqXfirstconstruct} && X \defas\sum_{j=1}^\infty
\zeta_j \delta_{\gamma_j},
\end{eqnarray}
where the random variables $\zeta_1, \zeta_2, \ldots$ are conditionally
independent given $B$\vspace*{-2pt} and
%
%e5 #&#
\begin{eqnarray}
\label{eqfirstBNBP} && \zeta_j \vert \basicbb_j \sim
\nbdist(r, \basicbb_j), \qquad j\in\mathbb{N},
\end{eqnarray}
for some parameter $r>0$.
Here, $\nbdist(r,p)$ denotes the negative binomial distribution with
parameters $r>0$, $p\in(0,1]$, whose probability mass function
(p.m.f.)\vspace*{-2pt} is
%
%e6 #&#
\begin{eqnarray}
\nbdist(z; r,p)\defas\frac{ (r)_z }{ z! } p^z (1-p)^r,
\qquad z\in\mathbb{Z}_+,
\end{eqnarray}
where $(a)_n \defas a(a+1)\cdots(a+n-1)$ with $(a)_0 \defas1$ is the
$n$th rising factorial.
Note that, conditioned on $B$, the point process $X$ is the (fixed
component) of a \emph{negative binomial process} \cite{BMPJpre,ZHDC12}.
Unconditionally, $X$ is the ordinary component of a beta negative
binomial process, which we formally define in~Section~\ref{secprelims}.

Conditioned on $B$, construct a sequence of point processes $\process X
{\mathbb{N}}$ that are i.i.d. copies of $X$.
In this case, $\process X {\mathbb{N}}$ is an exchangeable sequence of beta
negative binomial processes, and our primary goal is to characterize
the (unconditional) distribution of the sequence.
This task is non-trivial because the construction of the point process
$X$ in equation (\ref{eqXfirstconstruct}) is not \emph{finitary} in the sense
that no finite subset of the atoms of $B$ determines $X$ with
probability one.
In the case of conditionally-i.i.d. Bernoulli processes, the
unconditional distributions of the measures remain in the class of
Bernoulli processes, and so a finitary construction is
straightforwardly obtained with Poisson (point) processes.
Then the distribution of the sequence, which  Thibaux and Jordan \cite{TJ2007} showed is
characterized by the IBP, may be derived immediately from the conjugacy
between the classes of beta and Bernoulli processes \cite{Kim1999,Hjort1990,TJ2007}.
While conjugacy also holds between the classes of beta and negative
binomial processes \cite{BMPJpre,ZHDC12}, the unconditional law of
the point process $X$ is no longer that of a negative binomial process;
instead, it is the law of a beta negative binomial process.

Existing constructions for beta negative binomial processes truncate
the number of atoms in the underlying beta process and typically use
slice sampling to remove the error introduced by this approximation
asymptotically \cite{BMPJpre,ZHDC12,TGG07,PZWGC2010}.
In this work, we instead provide a construction for the beta negative
binomial process directly, avoiding a representation of the underlying
beta process.
To this end, note that while the beta process $B$ has a countably
infinite number of atoms \mbox{a.s.}, it can be shown that $B$ is
still an \mbox{a.s.}
finite measure \cite{Hjort1990}.
It follows as an easy consequence
that the point process $X$ is \mbox{a.s.} finite as well and,
therefore, has
an \mbox{a.s.} finite number of atoms, which we represent with a
Poisson process.
The atomic masses are then characterized by the \emph{digamma
distribution}, introduced by Sibuya \cite{Sibuya1979}, which has p.m.f. (for
parameters $r, \theta> 0$) given by\vspace*{-3pt}
%
%e7 #&#
\begin{eqnarray}
\label{eqintromixedbetapmf} && \digammadist(z; r, \theta) \defas \frac{1}{\psi(r+\theta) - \psi(\theta)}
\frac{(r)_z}{(r+\theta)_z} z^{-1}, \qquad z \ge 1,
\end{eqnarray}
where $\psi(a) \defas\Gamma'(a)/\Gamma(a)$ denotes the digamma function.
In Section~\ref{secNBPdefn}, we prove the following:
%
%th1 #&#
\begin{thm} \label{resultintromixedbnb}
Let $Y$ be a Poisson process on $(\Omega,\mathcal{A})$ with finite intensity
%
%e8 #&#
\begin{eqnarray}
\mathrm{d}s \mapsto c \bigl[ \psi( c + r ) - \psi( c ) \bigr] \widetilde
{B}_0( \mathrm{d}s ),
\end{eqnarray}
that is, $Y = \sum_{k=1}^\kappa\delta_{\gamma_k}$
for a Poisson random\vspace*{1pt} variable $\kappa$ with mean $c [ \psi( c + r ) -
\psi( c ) ] \widetilde{B}_0(\Omega)$ and i.i.d. random variables
$\nprocess
\gamma{\mathbb{N}}k$, independent from $\kappa$, each with distribution
$\widetilde{B}_0/\widetilde{B}_0(\Omega)$.
Let $\nprocess\zeta{\mathbb{N}}k$ be an independent collection of
i.i.d.
$\digammadist( r, c )$ random variables.
Then
%
%e9 #&#
\begin{eqnarray}
\label{eqXsecondconstruct} && X \stackrel{d} {=}\sum_{k=1}^\kappa
\zeta_k \delta_{\gamma_k},
\end{eqnarray}
where $X$ is the beta negative binomial process defined in equation (\ref{eqXfirstconstruct}).
\end{thm}

With this construction and conjugacy (the relevant results are
reproduced in Section~\ref{seccombstructure}), characterizing the
distribution of $\process X {\mathbb{N}}$ is straightforward.
However, in applications we are only interested in the \emph{combinatorial structure} of the sequence $\process X {\mathbb{N}}$,
that is,
the pattern of sharing amongst the atoms while ignoring the locations
of the atoms themselves.
More precisely, for every $n\in\mathbb{N}$, let $\Histories_n \defas
\mathbb{Z}_+
^n \setminus\theset{0^n}$
be the set of all length-$n$ sequences of non-negative integers,
excluding the all-zero sequence.
Elements in $\Histories_n$ are called \emph{histories}, and can be
thought of as representations of non-empty multisets of
$[n] \defas\{1,\ldots,n\}$.
For every $h \in\Histories_n$, let $\historycount{h}$ be the number of
elements $s \in\Omega$ such that $X_j\theset s = h(j)$ for all $j
\le n$.
By the combinatorial structure of a finite subsequence $X_{[n]} \defas
(X_1, \ldots, X_n)$, we will mean the collection $\gprocessM{n}$ of
counts, which together can be understood as representations of
multisets of histories.
These counts are combinatorial in the following sense: Let $\phi\colon
(\Omega,\mathcal{A})\to(\Omega,\mathcal{A})$ be a Borel
automorphism on $(\Omega,\mathcal{A})$,
that is, a measurable permutation of $\Omega$ whose inverse is also
measurable, and define the transformed processes $X_j^\phi\defas X_j
\circ\phi^{-1}$, for every\vspace*{1pt} $j\le n$, where each atom $s$ is
repositioned to $\phi(s)$.
The collection $\gprocess{\historycount{h}}{h\in\Histories_n}$ is
invariant to this transformation, and it is in this sense that they
only capture the combinatorial structure.
In Section~\ref{seccombstructure}, we prove the following.

%th2 #&#
\begin{thm} \label{resultNBIBPPMF}
The probability mass function of $\gprocessM{n}$ is
\begin{eqnarray}
&&\mathbb{P}\{ \historycount{h} = m_h \colon h\in
\Histories_n \}
\nonumber
\\[-13pt]
\label{eqNBIBPPMF}
\\[-1pt]
\nonumber
&& \quad = \frac{ (c T)^{ \sum_{h\in\Histories_n} m_h } }{
\prod_{h\in\Histories_n} m_h !} \exp \bigl( - c T\bigl[ \psi( c + n r ) -
\psi( c ) \bigr] \bigr) \prod_{h\in\Histories_n} \Biggl[
\frac{ \Gamma(\historysum{h}) \Gamma(c+nr) }{ \Gamma
(c+nr+\historysum
{h}) } \prod_{j=1}^n
\frac{ (r)_{h(j)} }{ h(j) ! } \Biggr] ^{m_h},
\nonumber
\end{eqnarray}
where $\historysum{h} \defas\sum_{j\le n} h(j)$, for every $h\in
\Histories_n$,
and $T\defas\widetilde{B}_0(\Omega) > 0$.
\end{thm}
%
%%%%%%%%%%%%%%%%%%%%%%%%%%%%%%%%%%%%%%%%%%%%%%%%%

As one would expect, equation (\ref{eqNBIBPPMF}) is reminiscent of the p.m.f.
for the IBP, and indeed the collection $\gprocess{\historycount
{h}}{h\in
\Histories_n}$ is characterized by what we call the \emph{negative
binomial Indian buffet process}, or NB-IBP.
Let $\bnbdist(r,\alpha,\beta)$ denote the \emph{beta negative
binomial distribution (with parameters $r, \alpha, \beta>0$)}, that
is, we write
$Z\sim\bnbdist(r,\alpha,\beta)$ if there exists a beta random
variable $p \sim\betadist(\alpha,\beta)$
such that ${Z \vert p \sim\nbdist(r, p)}$.
In the NB-IBP, a sequence of customers enters an Indian buffet restaurant:
\begin{itemize}
\item The first customer
\begin{itemize}
\item selects
$\poisson ( c \gamma[ \psi(c+r) - \psi(c) ]  )$
distinct dishes, taking $\digammadist(r,c)$ servings of each dish,
independently.
\end{itemize}

\item For $n\ge1$, the ($n+1$)st customer
\begin{itemize}
\item takes
$\bnbdist( r, S_{n,k}, c+nr )$
servings of each previously sampled dish $k$;
where $S_{n,k}$ is the total number of servings taken of dish $k$ by
the first $n$ customers;

\item selects
$\poisson ( c\gamma[ \psi(c+(n+1)r) - \psi(c+nr) ]  )$
new dishes to taste, taking  $\digammadist( r, c+nr )$ servings of each
dish, independently.
\end{itemize}
\end{itemize}

The interpretation here is that, for every $h \in\Histories_n$, the
count $\historycount{h}$ is the number of dishes $k$ such that, for
every $j\le n$, customer $j$ took $h(j)$ servings of dish $k$.
Then the sum $\historysum{h}$ in equation~(\ref{eqNBIBPPMF}) is the total
number of servings taken of dish $k$ by the first $n$ customers.
Because the NB-IBP is the combinatorial structure of a conditionally
i.i.d. process, its distribution, given in  Theorem~\ref{resultNBIBPPMF},
must be invariant to every permutation of the customers. We can state
this property formally as follows.
%
%th3 #&#
\begin{thm}[(Exchangeability)] \label{resultNBIBPexch}
Let $\pi$ be a permutation of $[n] \defas\{1,\ldots,n\}$, and, for
$h\in\Histories_n$,
note that the composition $h \circ\pi\in\Histories_n$ is given by
$(h \circ\pi)(j) = h(\pi(j))$, for every $j \le n$.
Then
%
%e10 #&#
\begin{eqnarray}
&& \gprocess{\historycount{h}} {h\in\Histories_n} \stackrel{d} {=}
\gprocess{\historycount{h \circ\pi}} {h\in\Histories_n}.
\end{eqnarray}
\end{thm}

The exchangeability of the combinatorial structure and its p.m.f. in
equation (\ref{eqNBIBPPMF}) allows us to develop Gibbs sampling techniques
analogous to those originally developed for the IBP \cite{GGS2007,meeds2007modeling}.
In particular, because the NB-IBP avoids a representation of the beta
process underlying the exchangeable sequence $\process X {\mathbb{N}}$, these
posterior inference algorithms do not require numerical integration
over representations of the beta process.
We discuss some of these techniques in Section~\ref{secinference}.

%s2 #&#
\section{Preliminaries}
\label{secprelims}

Here, we review \emph{completely random measures} and formally define
the negative binomial and beta negative binomial processes.
We provide characterizations via Laplace functionals and conclude the
section with a discussion of related work.

%%%%%%%%%%%%%%%%%%%%%%%%%%%%%
%s2.1 #&#
\subsection{Completely random measures}
\label{secCRM}
%%%%%%%%%%%%%%%%%%%%%%%%%%%%%

%
Let $\Mspace$ denote the space of $\sigma$-finite measures on
$(\Omega,\mathcal{A})$ equipped with the $\sigma$-algebra generated
by the
projection maps $\mu\mapsto\mu(A)$ for all $A\in\mathcal{A}$.
A random measure $\xi$ on $(\Omega,\mathcal{A})$ is a random element in
$\Mspace$, and we say that $\xi$ is \emph{completely random} or
\emph{has independent increments} when, for every finite collection of
disjoint,
measurable sets $A_1,\ldots,A_n \in\mathcal{A}$, the random variables
$\xi(A_1),\ldots,\xi(A_n)$
are independent.
Here, we briefly review completely random measures; for a thorough
treatment, the reader should consult Kallenberg \cite{Kallenberg2002}, Chapter~12,  or the classic text by
Kingman \cite{Kingman1967}.
Every completely random measure $\xi$ can be written as a sum of
three independent parts\vspace*{-3pt}
%
%e11 #&#
\begin{eqnarray}
\label{eqgeneralCRM} && \xi= \bar\xi+ \sum_{s\in\Atoms}
\vartheta_s \delta_s + \sum
_{(s,p) \in\eta} p \delta_s \qquad\mbox{a.s.},
\end{eqnarray}
called the \emph{diffuse}, \emph{fixed}, and \emph{ordinary}
components, respectively, where:

\begin{longlist}[3.]
\item[1.] $\bar\xi$ is a non-random, non-atomic measure;

\item[2.] $\Atoms\subseteq\Omega$ is a non-random countable set whose
elements are referred to as the
\emph{fixed atoms} and whose masses $\vartheta_1, \vartheta_2,
\ldots$
are independent random variables in $\mathbb{R}_+$ (the non-negative
real numbers);

\item[3.] $\eta$ is a Poisson process on $\Omega\times(0,\infty)$ whose
intensity $\mathbb{E}\eta$ is $\sigma$-finite and has diffuse projections
onto $\Omega$, that is, the measure $(\mathbb{E}\eta) (\cdot
\times
(0,\infty))$ on $\Omega$ is non-atomic.
\end{longlist}
In this article, we will only study purely-atomic completely random
measures, which therefore have no diffuse component.
It follows that we may characterize the law of $\xi$ by (1) the
distributions of the atomic masses in the fixed component, and (2) the
intensity of the Poisson process underlying the ordinary\vspace*{-3pt} component.

%s2.2 #&#
\subsection{Definitions}

By a \emph{base measure} on $(\Omega,\mathcal{A})$, we mean a
$\sigma$-finite measure $B$ on $(\Omega,\mathcal{A})$ such that $B \{
s\} \le1$
for all $s \in\Omega$.
For the remainder of the article, fix a base measure $B_0$. We may\vspace*{-3pt} write
%
%e12 #&#
\begin{eqnarray}
\label{eqbasemeas} && B_0= \widetilde{B}_0+ \sum
_{s\in\Atoms} \bar{b}_s \delta_s
\end{eqnarray}
for some
non-atomic measure $\widetilde{B}_0$;
a countable set $\Atoms\subseteq\Omega$;
and constants $\bar{b}_1, \bar{b}_2, \ldots$ in $(0,1]$.\footnote
{Note that we
have relaxed the condition on $\widetilde{B}_0$ (in the \hyperref[sec1]{Introduction})
to be
merely $\sigma$-finite.}
As discussed in the \hyperref[sec1]{Introduction}, a convenient model for random base
measures are \emph{beta processes}, a class of completely random
measures introduced by Hjort \cite{Hjort1990}.
For the remainder of the article, let $c \colon\Omega\to\mathbb{R}_+$
be a measurable function, which we call a \emph{concentration function}
(or \emph{parameter} when it is constant).

%
%de1 #&#
\begin{definition}[(Beta process)]
A random measure $B$ on $(\Omega,\mathcal{A})$ is a \emph{beta
process with concentration function $c$ and base measure $B_0$},
written $B\sim
\mathrm{BP}_{\mathcal{L}}( c, B_0)$, when it is purely atomic and
completely random,
with a fixed\vspace*{-3pt} component
%
%e13 #&#
\begin{eqnarray}
\label{eqBPfixed} && \sum_{s\in\Atoms} \vartheta_s
\delta_s, \qquad \vartheta_s \stackrel{\mathrm{ind}} {\sim}
\betadist\bigl(c(s) \bar{b}_s, c(s) (1 - \bar{b}_s)
\bigr),
\end{eqnarray}
and an ordinary component with intensity\vspace*{-3pt} measure
%
%e14 #&#
\begin{eqnarray}
\label{eqBPintensity} &&(\mathrm{d}s, \mathrm{d}p) \mapsto c(s)p^{-1} (1 -
p)^{c(s)-1}\, \mathrm{d}p \widetilde{B}_0(\mathrm{d}s).
\end{eqnarray}
\end{definition}

It is straightforward to show that a beta process is itself a base
measure with probability one.
This definition of the beta process generalizes the version given in
the introduction to a non-homogeneous process with a fixed component.
Likewise, we generalize our earlier definition of a \emph{negative
binomial process} to include an ordinary component.

%de2 #&#
\begin{definition}[(Negative binomial process)] \label{defnbp}
A point process $X$ on $(\Omega,\mathcal{A})$ is a \emph{negative
binomial process with parameter $r>0$ and base measure $B_0$}, written
$X \sim
\NBPLAW(r, B_0)$, when it is purely atomic and completely random,
with a fixed\vspace*{-2pt} component
%
%e15 #&#
\begin{eqnarray}
&& \sum_{s\in\Atoms} \vartheta_s
\delta_s, \qquad\vartheta_s \stackrel{\mathrm{ind}} {\sim}
\nbdist( r, \bar{b}_s ),
\end{eqnarray}
and an ordinary component with intensity measure\vspace*{-2pt}
%
%e16 #&#
\begin{eqnarray}
&& (\mathrm{d}s, \mathrm{d}p) \mapsto r \delta_1(\mathrm{d}p)
\widetilde{B}_0(\mathrm{d}s).
\end{eqnarray}
\end{definition}

The fixed component in this definition was given by
Broderick \textit{et al.} \cite{BMPJpre} and
Zhou \textit{et al.} \cite{ZHDC12} (and by Thibaux  \cite{ThibauxThesis} for the case \mbox{$r=1$}).
Here, we have additionally defined an ordinary component, following
intuitions from Roy \cite{Roy13CUP}.

The law of a random measure is completely characterized by its Laplace
functional, and this representation is often simpler to manipulate:
From Campbell's theorem, or a version of the L\'evy--Khinchin formula
for Borel spaces,
one can show that the Laplace functional of $X$ is\vspace*{-2pt}
%
%e17 #&#
\begin{eqnarray}\label{resultcfnbp}
&& f \mapsto \mathbb{E}\bigl[ e^{-X(f) } \bigr] = \exp \biggl[ - \int \bigl(
1 - e^{- f(s)} \bigr) r \widetilde{B}_0( \mathrm{d}s ) \biggr]
\prod_{s\in\Atoms} \biggl[ \frac{ 1- \bar{b}_s } { 1 - \bar{b}_s e^{-f(s)} }
\biggr]^r,
\end{eqnarray}
where $f$ ranges over non-negative measurable functions and
$X(f) \defas\int f(s) X(\mathrm{d}s)$.

Finally, we define \emph{beta negative binomial processes} via their
conditional law.

%de3 #&#
\begin{definition}[(Beta negative binomial process)] \label{defnBNBP}
A random measure $X$ on $(\Omega,\mathcal{A})$ is a \emph{beta
negative binomial process with parameter $r>0$, concentration function
$c$, and base measure $B_0$}, written\vspace*{-2pt}
\begin{eqnarray}
&& X \sim\BNBPLAW(r, c, B_0),
\nonumber
\end{eqnarray}
if there exists a beta process $B\sim\mathrm{BP}_{\mathcal{L}}( c, B_0)$
such\vspace*{-2pt} that
%
%e18 #&#
\begin{eqnarray}
&& X \vert B \sim\NBPLAW(r, B).
\end{eqnarray}
\end{definition}

This characterization was given by Broderick \textit{et al.} \cite{BMPJpre} and can be seen to
match a special case of the model in Zhou \textit{et al.} \cite{ZHDC12} (see the discussion\
of related work in Section~\ref{secrelatedworks}). It is straightforward to
show that a beta negative binomial process is also completely random,
and that its Laplace functional is given by
\begin{eqnarray}
\mathbb{E}\bigl[ e^{- X(f) } \bigr] & =&  \exp \biggl[ - \int \biggl[ 1 -
\biggl( \frac{ 1-p }{ 1-p e^{-f(s)} } \biggr)^r  \biggr] c(s) p^{-1}
(1-p)^{c(s)-1} \,\mathrm{d}p \widetilde{B}_0(\mathrm{d}s)
\biggr]
\nonumber
\\[-8pt]
\label{resultcfbnbp}
\\[-8pt]
\nonumber
&& {}\times \prod_{s\in\Atoms} \int \biggl(
\frac{ 1-p }{ 1- p e^{-f(s)} } \biggr)^r \betadist\bigl(p; c(s)
\bar{b}_s, c(s) (1-\bar{b}_s)\bigr) \,\mathrm{d}p,
\nonumber
\end{eqnarray}
for $f\colon\Omega\to\mathbb{R}_+$ measurable, where we note that the
factors in the product term take the form of the Laplace transform of
the beta negative binomial distribution.

%s2.3 #&#
\subsection{Related work}
\label{secrelatedworks}

The term ``negative binomial process'' has historically been reserved
for processes with negative binomial increments -- a class into which
the process we study here does not fal -- and these processes have been
long-studied in probability and statistics.
We direct the reader to Kozubowski and Podg{\'o}rski  \cite{Kozubowski2009}
for references.

One way to construct a process with negative binomial increments is to
rely upon the fact that a negative binomial distribution is a gamma
mixture of Poisson distributions. In particular, similarly to the
construction by Lo  \cite{lo1982bayesian}, consider a Cox process $X$
directed by a gamma process $G$ with finite non-atomic intensity. So
constructed, $X$ has independent increments with negative binomial
distributions.
Like the beta process (with a finite intensity underlying its ordinary
component), the gamma process has, with probability one, a countably
infinite number of atoms but a finite total mass, and so the Cox
process $X$ is \mbox{a.s.} finite as well. Despite similarities, a comparison
of Laplace functionals shows that the law of $X$ is not that of a beta
negative binomial process.
Using an approach directly analogous to the derivation of the IBP in
\cite{GG06}, Titsias \cite{Titsias2007} characterizes the combinatorial
structure of a sequence of point processes that, conditioned on $G$,
are independent and identically distributed to the Cox process $X$. See
Section~\ref{seccombstructure} for comments.
This was the first count analogue of the IBP; the possibility of a
count analogue arising from beta negative binomial processes was first
raised by Zhou \textit{et al.} \cite{ZHDC12},
who described the distribution of the number of new dishes sampled by
each customer.
Recent work by  Zhou, Madrid and Scott \cite{ZMS2014}, independent of our own and proceeding
along different lines, describes a combinatorial process related to the
NB-IBP (following a re-scaling of the beta process intensity).

Finally, we note that another negative binomial process without
negative binomial increments was defined on Euclidean space by Barndorff-Nielsen and Yeo  \cite{BY1969}
and extended to general spaces by Gr{\'e}goire \cite{Gregoire1984}
and
Wolpert and Ickstadt \cite{WI1998Poisson}.
These measures are generally Cox processes on $(\Omega,\mathcal{A})$
directed by
random measures of the form
\begin{eqnarray}
&& \mathrm{d}s \mapsto\int_{\mathbb{R}_+} \nu(t,\mathrm{d}s) G (
\mathrm{d}t),
\nonumber
\end{eqnarray}
where $G$ is again a gamma process, this time on $\mathbb{R}_+$, and
$\nu$
is a probability kernel from $\Omega$ to $\mathbb{R}_+$, for example, the
Gaussian kernel.

%s3 #&#
\section{Constructing beta negative binomial processes}
\label{secNBPdefn}

Before providing a finitary construction for the beta negative binomial
process, we make a few remarks on the digamma distribution.
For the remainder of the article, define $\harmonicsum{r,\theta}
\defas
\psi(\theta+r)-\psi(\theta)$ for some $r, \theta> 0$.
Following a representation by Sibuya \cite{Sibuya1979}, we may relate the
digamma and beta negative binomial distributions as follows:
Let $Z \sim\digammadist( r, \theta)$ and define $W\defas Z-1$, the
latter of which has p.m.f.
%
%e19 #&#
\begin{eqnarray}\label{eqdigammarewrite}
&& \mathbb{P}\{W=w\} = (\theta\harmonicsum{r,\theta})^{-1}
\frac{ w+r }{ w+1 } \bnbdist(w;r,1,\theta), \qquad w\in\mathbb{Z}_+.
\end{eqnarray}
Deriving\vspace*{1pt} the Laplace transform of the law of $W$ is straightforward,
and because $\mathbb{E}[ e^{-tW} ] = e^t \mathbb{E}[ e^{-tZ} ]$, one
may verify
that the Laplace transform of the digamma distribution is given by
%
%e20 #&#
\begin{eqnarray}\label{eqlaplacemixedbnb}
&& \Psi_{r,\theta}(t) \defas \mathbb{E}\bigl[ e^{-tZ} \bigr] = 1 -
\harmonicsum{r,\theta}^{-1} \int \biggl[ 1- \biggl( \frac{1-p}{1-p  e^{-t}}
\biggr)^r \biggr]  p^{-1} (1-p)^{\theta-1} \,\mathrm{d}p.
\end{eqnarray}

The form of equation (\ref{eqdigammarewrite}) suggests the following rejection
sampler, which was first proposed
by Devroye \cite{devroye1992random}, Proposition~2, Remark~1:
Let $r>0$ and let $\process U {\mathbb{N}}$ be an i.i.d. sequence of uniformly
distributed random numbers.
Let
\begin{eqnarray}
&& \process Y {\mathbb{N}}\stackrel{\mathrm{i.i.d.}} {\sim}\bnbdist(r, 1, \theta),
\nonumber
\end{eqnarray}
and define $\eta\defas\inf \theset{ n \in\mathbb{N}:\max
\theset{r,1} \cdot U_n < \frac{Y_i+r} {Y_i+1} }$. Then
\begin{eqnarray}
&& Y_\eta+ 1 \sim\digammadist(r, \theta),
\nonumber
\end{eqnarray}
and
\begin{eqnarray}
&& \mathbb{E}\eta= \frac{ \max\{r,1\} }{ \theta[ \psi(r + \theta) -
\psi
(\theta
) ] }; \qquad \mathbb{E}\eta< \max
\theset{r,r^{-1}}.
\nonumber
\end{eqnarray}

With digamma random variables, we provide a finitary construction for
the beta negative binomial process.
The following result generalizes the statement given by Theorem~\ref
{resultintromixedbnb} (in the \hyperref[sec1]{Introduction}) to a non-homogeneous
process, which also has a fixed component.

%th4 #&#
\begin{thm} \label{resultmixedbnb}
Let $r>0$, and let $\vartheta\defas\gprocess{\vartheta_s}{s\in
\Atoms
}$ be a collection of independent random variables with
%
%e21 #&#
\begin{eqnarray}
&& \vartheta_s \sim\bnbdist\bigl( r, c(s) \bar{b}_s,
c(s) ( 1- \bar{b}_s ) \bigr), \qquad s\in\Atoms.
\end{eqnarray}
Let $Y$ be a Poisson process on $(\Omega,\mathcal{A})$, independent from
$\vartheta$, with (finite) intensity
%
%e22 #&#
\begin{eqnarray}
&& \mathrm{d}s \mapsto c(s) \bigl[ \psi\bigl( c(s) + r \bigr) - \psi\bigl( c(s)
\bigr) \bigr] \widetilde{B}_0( \mathrm{d}s).
\end{eqnarray}
Write $Y = \sum_{k=1}^\kappa\delta_{\gamma_k}$ for some random element
$\kappa$ in $\mathbb{Z}_+$ and \mbox{a.s.} unique random elements
$\gamma_1,
\gamma
_2, \ldots$ in $\Omega$, and put $\mathcal{F}\defas\sigma(\kappa,
\gamma
_1, \gamma_2, \ldots)$.
Let $\gprocess{\zeta_j}{j\in{\mathbb{N}}}$ be a collection of random variables
that are independent from $\vartheta$ and are conditionally independent
given $\mathcal{F}$, and let
%
%e23 #&#
\begin{eqnarray}\label{eqordvar}
&&\zeta_j \vert \mathcal{F} \sim\digammadist\bigl( r, c(
\gamma_j) \bigr), \qquad j\in{\mathbb{N}}.
\end{eqnarray}
Then
%
%e24 #&#
\begin{eqnarray}
&& X = \sum_{s\in\Atoms} \vartheta_s
\delta_s + \sum_{j=1}^\kappa
\zeta_j \delta_{\gamma_j} \sim\BNBPLAW( r, c, B_0).
\end{eqnarray}
\end{thm}

\begin{pf}
We have
%
%e25 #&#
\begin{eqnarray}
&& \mathbb{E}^\mathcal{F}\bigl[ e^{-X(f)} \bigr] = \prod
_{s\in\Atoms} \mathbb{E}\bigl[ e^{- \vartheta_s f(s)} \bigr] \times \prod
_{j=1}^\kappa \mathbb{E}^\mathcal{F}
\bigl[ e^{- \zeta_j f(\gamma_j) } \bigr],
\end{eqnarray}
for every $f\colon\Omega\to\mathbb{R}_+$ measurable.
For $s \in\Omega$,
write
$
g(s) = \Psi_{r,c(s)}  (f(s)  )
$
for the Laplace transform of the digamma distribution evaluated at $f(s)$,
where $\Psi_{r,\theta}(t)$ is given by equation (\ref{eqlaplacemixedbnb}).
We may then write
%
%e26 #&#
\begin{eqnarray}
&& \prod_{j=1}^\kappa \mathbb{E}^\mathcal{F}
\bigl[ e^{-\zeta_j f(\gamma_j)} \bigr] = \prod_{j=1}^\kappa
g( \gamma_j ).
\end{eqnarray}
Then by the chain rule of conditional expectation, complete randomness,
and Campbell's theorem,
%
%e27 #&#
%e28 #&#
\begin{eqnarray}
\mathbb{E}\bigl[ e^{-X(f)} \bigr] &=& \prod_{s\in\Atoms}
\mathbb{E}\bigl[ e^{-\vartheta_s f(s)} \bigr] \times \exp \biggl[ - \int
_\Omega\bigl(1-g(s)\bigr) c(s) \harmonicsum{r,c(s)} \widetilde
{B}_0(\mathrm{d}s) \biggr]
\\
& =& \prod_{s\in\Atoms} \biggl[ \int \biggl(
\frac{ 1-p }{ 1- p e^{-f(s)} } \biggr)^r \betadist\bigl(p; c(s)
\bar{b}_s, c(s) (1-\bar{b}_s)\bigr) \,\mathrm{d}p \biggr]
\\
&&{}\times \exp \biggl[ - \int_{(0,1] \times\Omega} \biggl[ 1 - \biggl(
\frac{ 1-p }{ 1-p e^{-f(s)} } \biggr)^r \biggr] c(s) p^{-1}
(1-p)^{c(s)-1} \,\mathrm{d}p \widetilde{B}_0(\mathrm{d}s)
\biggr],
\nonumber
\end{eqnarray}
which is the desired form of the Laplace functional.
\end{pf}

A finitary construction for conditionally-i.i.d. sequences of negative
binomial processes with a common beta process base measure now follows
from known conjugacy results.
In particular, for every $n\in{\mathbb{N}}$, let $\Xseqn{n} \defas(X_1,
\ldots
, X_n)$.
The following theorem characterizes the conjugacy between the (classes
of) beta and negative binomial processes and follows from repeated
application of the results by Kim \cite{Kim1999}, Theorem~3.3  or
Hjort \cite{Hjort1990}, Corollary~4.1.
This result, which is tailored to our needs, is similar to those
already given by  Broderick \textit{et el.} \cite{BMPJpre}
and  Zhou \textit{et al.} \cite{ZHDC12}, and generalizes
the result given by  Thibaux \cite{ThibauxThesis} for the case $r=1$.
%
%th5 #&#
\begin{thm}[(Hjort \cite{Hjort1990}, Zhou \textit{et al.} \cite{ZHDC12})]
\label{resultbnbpconjugacy}
Let $B \sim\mathrm{BP}_{\mathcal{L}}(c,B_0)$ and, conditioned on
$B$, let $\process X
{\mathbb{N}}$ be a sequence of i.i.d. negative binomial processes with
parameter $r>0$ and base measure $B$.
Then for every $n\in{\mathbb{N}}$,
%
%e29 #&#
\begin{eqnarray}
&& B \vert X_{[n]} \sim\mathrm{BP}_{\mathcal{L}} \biggl( c_n,
\frac{c}{c_n} B_0+ \frac{1}{c_n} S_n \biggr),
\end{eqnarray}
where $S_n \defas\sum_{i=1}^n X_i$ and $c_n(s) \defas c(s) +
S_n\theset s + nr$, for $s \in\Omega$.
\end{thm}

%re1 #&#
\begin{remark} \label{resultbnbppredictive}
It follows immediately that, for every $n\in{\mathbb{N}}$, the law of
$X_{n+1}$ conditioned on $X_1,\ldots,X_n$ is given by
%
%e30 #&#
\begin{eqnarray}
\label{eqbnbppredictive}
X_{n+1} \vert X_{[n]} \sim\BNBPLAW \biggl( r,
c_n, \frac{c}{c_n} B_0+ \frac{1}{c_n}
S_n \biggr).
\end{eqnarray}
\end{remark}

We may therefore construct this exchangeable sequence of beta negative
binomial processes with Theorem~\ref{resultmixedbnb}.

%%%%%%%%%%%%%%%%%%%%%%%%%%%%%%%
%s4 #&#
\section{Combinatorial structure}
\label{seccombstructure}
%%%%%%%%%%%%%%%%%%%%%%%%%%%%%%%

%
We now characterize the combinatorial structure of the exchangeable
sequence $X_{[n]}$ in the case when $c>0$ is constant and $B_0(=
\widetilde{B}_0
)$ is non-atomic.
In order to make this precise, we introduce a quotient of the space of
sequences of integer-valued measures.
Let $n\in{\mathbb{N}}$ and for any pair $U \defas(U_1, \ldots, U_n)$
and $V
\defas(V_1, \ldots, V_n)$ of (finite) sequences of integer-valued
measures, write $U \sim V$ when there exists a Borel automorphism $\phi
$ on $(\Omega,\mathcal{A})$ satisfying $U_j = V_j \circ\phi^{-1}$
for every
$j\le n$. It is easy to verify that $\sim$ is an equivalence relation.
Let $\CombStruct{U}$ denote the equivalence class containing $U$.
The quotient space induced by $\sim$ is itself a Borel space, and can
be related to the Borel space of sequences of $\mathbb{Z}_+$-valued measures
by coarsening the $\sigma$-algebra to that generated by the functionals
%
%e31 #&#
\begin{eqnarray}
&& \mathcal M_h (U_1, \ldots, U_n ) \defas
\cardsharp\bigl\{ s\in\Omega:\forall j\le n,  U_j\theset s = h(j)
\bigr\}, \qquad  h\in\Histories_n, j\le n,
\end{eqnarray}
where $\cardsharp A$ denotes the cardinality of $A$, and $\Histories_n
\defas\mathbb{Z}_+^n \setminus\theset{0^n}$ is the space of histories
defined in the \hyperref[sec1]{Introduction}.
The collection $\gprocessM{n}$ of multiplicities (of histories)
corresponding to $\Xseqn{n}$, also defined in the \hyperref[sec1]{Introduction}, then satisfies
$\historycount{h} = \mathcal M_h ( \Xseqn{n} )$ for every $h\in
\Histories_{n}$.
The collection $\gprocessM{n}$ thus identifies a point in the quotient
space induced by $\sim$.
Our aim is to characterize the distribution of $\gprocessM{n}$, for
every $n \in{\mathbb{N}}$.

Let $\oh\in\Histories_n$, and define $\Histories_{n+1}^{(\oh)}
\defas\{ h\in\Histories_{n+1} :\forall j \le n,  h(j) = \oh(j)
\}
$ to be the collection of histories in $\Histories_{n+1}$ that agree
with $\oh$ on the first $n$ entries.
Then note that
%
%e32 #&#
\begin{eqnarray}\label{eqDetermineCounts}
&& \historycount{\oh} = \sum_{ h\in\Histories_{n+1}^{(\oh)} } \historycount{h},
\qquad\oh\in\Histories_n,
\end{eqnarray}
that is, the multiplicities $\gprocess{\historycount{h}}{h\in
\Histories
_{n+1}}$ at stage $n+1$ completely determine the multiplicities
$\gprocess{\historycount{\oh}}{\oh\in\Histories_n}$ at all
earlier stages.
It follows that
%
%e33 #&#
\begin{eqnarray}
\mathbb{P}\{ \historycount{h} = m_h \colon h\in
\Histories_{n+1} \} &= & \mathbb{P}\{ \historycount{\oh} = m_{\oh}
\colon\oh\in \Histories_n \}
\nonumber
\\[-8pt]
\label{eqJointCombStruct}
\\[-8pt]
\nonumber
&&{}\times \mathbb{P}\{ \historycount{h} = m_h \colon h
\in\Histories_{n+1} \vert \historycount{\oh} = m_{\oh
} \colon\oh
\in\Histories_n \},
\end{eqnarray}
where $m_{\oh} = \sum_{ h\in\Histories_{n+1}^{(\oh)} } m_h$ for
$\oh
\in\Histories_n$.
The structure of equation (\ref{eqJointCombStruct})
suggests an
inductive proof for Theorem~\ref{resultNBIBPPMF}.

%%%%%%%%%%%%%%%%%%%%%%%%%%%%%%%
%s4.1 #&#
\subsection{The law of \texorpdfstring{$\historycount{h}$}{$M_m$} for \texorpdfstring{$h\in\Histories_1$}{$h in H_1$}}
%%%%%%%%%%%%%%%%%%%%%%%%%%%%%%%

Note that $\Histories_1$ is isomorphic to ${\mathbb{N}}$ and that the
collection $\gprocess{\historycount{h}}{h\in\Histories_1}$ counts the
number of atoms of each positive integer mass.
It follows from Theorem~\ref{resultintromixedbnb} and a transfer argument
\cite{Kallenberg2002}, Propositions~6.10, 6.11 and 6.13,  that there exists:
\begin{longlist}[3.]
\item[1.] a Poisson random variable $\kappa$ with mean $c T
\harmonicsum{r, c} $, where $T\defas\widetilde{B}_0(\Omega) <
\infty$;
\item[2.] an i.i.d. collection of a.s. unique random elements $\gamma
_{1}, \gamma_{2}, \ldots$ in $\Omega$;
\item[3.] an i.i.d. collection $\gprocess{\zeta_{j}}{j\in{\mathbb{N}}}$ of
$\digammadist(r,c)$ random variables;
\end{longlist}
all mutually independent,
such that
\begin{eqnarray}
X_1 &=& \sum_{j=1}^{\kappa}
\zeta_{j} \delta_{\gamma_{j}} \qquad\mbox{a.s.}
\nonumber
\end{eqnarray}
It follows that
%
%e34 #&#
\begin{eqnarray}
\historycount{h} &=& \cardsharp\bigl\{ j\le\kappa\colon\zeta_{j} = h(1)
\bigr\} \qquad\mbox{a.s.}, \mbox{ for } h\in\Histories_1,
\end{eqnarray}
and $\kappa= \sum_{h\in\Histories_1} \historycount{h}$ a.s.
Therefore,
%
%e35 #&#
\begin{eqnarray}
&& \mathbb{P}\{ \historycount{h} = m_h \colon h\in
\Histories_1 \}
\nonumber
\\[-8pt]
\\[-8pt]
\nonumber
&&\quad= \mathbb{P}\biggl\{ \kappa= \sum_{h\in\Histories_1}
m_h\biggr\} \mathbb{P}\biggl\{ \historycount{h} = m_h
\colon h\in\Histories_1 \Big\vert \kappa= {\sum
_{h\in\Histories_1} m_h} \biggr\}.
\end{eqnarray}
Because $\zeta_1,\zeta_2,\ldots$ are i.i.d.,
the collection $\gprocessM{1}$ has a multinomial distribution
conditioned on its sum $\kappa$. Namely,
$M_h$ counts the number of times, in $\kappa$ independent trials, that
the multiplicity $h(1)$ arises from a $\digammadist(r,c)$ distribution.
In particular,
%
%e36 #&#
\begin{eqnarray}
&&\mathbb{P}\biggl\{ \historycount{h} = m_h \colon h\in
\Histories_1 \Big\vert \kappa= { \sum
_{h \in\Histories_1} m_h } \biggr\}
\nonumber
\\[-8pt]
\\[-8pt]
\nonumber
&&\quad = \frac{(\sum_{h\in\Histories_1} m_h)! }{ \prod_{h\in\Histories
_1} (
m_h ! ) } \prod_{h \in\Histories_1} \bigl[
\digammadist\bigl( h(1); r, c \bigr)^{m_h} \bigr]. %
\end{eqnarray}
It follows that
%
%e37 #&#
\begin{eqnarray}
&&\mathbb{P}\{ \historycount{h} = m_h \colon h\in
\Histories_1 \}
\nonumber
\\[-8pt]
\label{eqdistX1}
\\[-8pt]
\nonumber
&&\quad = \frac{ (c T\harmonicsum{r,c} )^{\sum_{h\in\Histories_1}
m_h } }{ { \prod_{h\in\Histories_1} ( m_h ! ) } } \exp( - c T\harmonicsum{r,c} ) \prod
_{h\in\Histories_1} \bigl[ \digammadist\bigl( h(1) ; r, c
\bigr)^{m_h} \bigr].
\end{eqnarray}

%%%%%%%%%%%%%%%%%%%%%%%%%%%%%%%
%s4.2 #&#
\subsection{The conditional law of $\historycount{h}$ for \texorpdfstring{$h \in\Histories_{n+1}$}{$h in\Histories_{n+1}$}}
%%%%%%%%%%%%%%%%%%%%%%%%%%%%%%%

Let $S_n \defas\sum_{j=1}^n X_j$.
Recall that $s(\oh) \defas\sum_{j\le n} \oh(j)$ for $\oh\in
\Histories_{n}$.
We may write
%
%e38 #&#
\begin{eqnarray}
S_n &=& \sum_{\oh\in\Histories_n} \sum
_{j=1}^{\historycount{\oh}} s(\oh ) \delta_{\omega_{\oh,j}},
\end{eqnarray}
for some collection $\omega\defas
\gprocess{\omega_{\oh,j}} {\oh\in\Histories_{n}, j \in{\mathbb{N}}}$
of \mbox{a.s.} distinct random elements in $\Omega$.
It follows from Remark~\ref{resultbnbppredictive},
Theorem~\ref{resultintromixedbnb}, and a transfer argument that there exists:
\begin{longlist}[4.]

\item[1.] a Poisson random variable $\kappa$ with mean $c T
\harmonicsum{r,c + nr} $;

\item[2.] an i.i.d. collection of \mbox{a.s.} unique random elements
$\gamma_1,
\gamma_2, \ldots$ in $\Omega$, \mbox{a.s.} distinct also from
$\omega$;

\item[3.] an i.i.d. collection $\gprocess{\zeta_{j}}{j\in{\mathbb{N}}}$ of
$\digammadist(r, c+nr)$ random variables;

\item[4.] for each $\oh\in\Histories_n$,
an i.i.d. collection $\gprocess{\vartheta_{\oh,j}}{j\in\mathbb
{N}}$ of random
variables satisfying
\begin{eqnarray}
&&  \vartheta_{\oh,j} \sim\bnbdist\bigl( r, s(\oh), c+nr \bigr) \qquad
\mbox{ for $j \in{\mathbb{N}}$;}
\nonumber
\end{eqnarray}
\end{longlist}
all mutually independent and independent of $X_{[n]}$,
such that
%
%e39 #&#
\begin{eqnarray} \label{eqcombgen}
X_{n+1} &=& \sum_{\oh\in\Histories_n} \sum
_{j=1}^{\historycount{\oh}} \vartheta_{\oh,j}
\delta_{\omega
_{\oh,j}} + \sum_{j=1}^{\kappa}
\zeta_{j} \delta_{\gamma_{j}} \qquad\mbox{a.s.}
\end{eqnarray}
Conditioned on $\Xseqn{n}$, the first and second terms on the right-hand side correspond to the fixed and ordinary components of $X_{n+1}$,
respectively.
Let
%
%e40 #&#
\begin{eqnarray}
&& \extendedspace_{n+1} \defas\bigl\{ h\in\Histories_{n+1} \colon
h(j) = 0, j\le n \bigr\}
\end{eqnarray}
be the set of histories $h$
for which $h(n+1)$ is the first non-zero element.
Then, with probability one,
%
%e41 #&#
\begin{eqnarray}
&& \historycount{h} = \cardsharp \bigl\{ j\le\kappa\colon\zeta_{j} =
h(n+1) \bigr\} \qquad\mbox{for $h\in\extendedspace_{n+1}$},
\end{eqnarray}
and
%
%e42 #&#
\begin{eqnarray}
\historycount{h} &=& \cardsharp \bigl\{ j \le\historycount{\oh} \colon
\vartheta_{\oh,j} = h(n+1) \bigr\} \qquad\mbox{for $\oh\in
\Histories_n$ and $h\in\extendedhist{n+1} {\oh}$}.
\end{eqnarray}
By the stated independence of the variables above, we have
%
%e43 #&#
\begin{eqnarray}
&&\mathbb{P}\{ \historycount{h} = m_h \colon h\in
\Histories_{n+1} \vert \historycount{\oh} = m_{\oh} \colon\oh\in
\Histories_n \}
\nonumber
\\[-8pt]
\label{eqCombStructGen}
\\[-8pt]
\nonumber
&&\quad = \mathbb{P}\bigl\{ \historycount{h} = m_h \colon h \in
\extendedspace _{n+1} \bigr\} \prod_{\oh\in\Histories_n}
\mathbb{P}\bigl\{ \historycount{h} = m_h \colon h\in
\Histories_{n+1}^{(\oh)} \vert \historycount{\oh} = m_{\oh}
\bigr\}.
\end{eqnarray}
Let $\existingspace_{n+1} \defas\bigcup_{\oh\in\Histories_{n}}
\extendedhist{n+1} {\oh}$.
For every $\oh\in\Histories_n$,
the random variables $\vartheta_{\oh,1},\vartheta_{\oh,2},\ldots$ are
i.i.d., and therefore, conditioned on
$\historycount{\oh}$,
the collection $\gprocess{\historycount{h}}{ h\in\Histories
_{n+1}^{(\oh
)} }$
has a multinomial distribution.
In particular, the product term in equation (\ref{eqCombStructGen}) is given by
%
%e44 #&#
\begin{eqnarray*}
&& \prod_{\oh\in\Histories_n} \mathbb{P}\bigl\{ \historycount{h} =
m_h \colon h\in \Histories_{n+1}^{(\oh)} \vert \historycount{\oh} = m_{\oh} \bigr\}
\\
&&\quad = \frac{ \prod_{\oh\in\Histories_n} (m_{\oh} !) }{ \prod_{h\in
\existingspace_{n+1}} (m_h !) } \prod_{h\in\existingspace_{n+1} } \bigl[
\bnbdist\bigl( h(n+1) ; r, \historysum{h} - h(n+1), c + nr \bigr)^{m_h}
\bigr].
\end{eqnarray*}
The p.m.f. of the beta negative binomial distribution is given by
%
%e45 #&#
\begin{eqnarray}
&&\bnbdist(z ; r, \alpha, \beta) = \frac{(r)_z }{z } \frac{ \betafn(z+\alpha, r+\beta) }{ \betafn(\alpha, \beta) },
\qquad z\in\mathbb{Z}_+,
\end{eqnarray}
for positive parameters $r, \alpha$, and $\beta$, where $\betafn
(\alpha
, \beta) \defas\Gamma(\alpha) \Gamma(\beta) / \Gamma(\alpha
+\beta)$
denotes the beta function.
We have that
$\kappa= \sum_{h\in\extendedspace_{n+1}} \historycount{h}$ \mbox{a.s.},
and therefore
%
%e46 #&#
\begin{eqnarray}
&& \mathbb{P}\bigl\{ \historycount{h} = m_h \colon h \in
\extendedspace _{n+1} \bigr\}
\nonumber
\\
&&\quad= \mathbb{P}\biggl\{ \kappa= { \sum_{h\in\extendedspace
_{n+1}}
m_h } \biggr\}
\\
&&\qquad{}\times \mathbb{P}\biggl\{ \historycount{h} =
m_h \colon h \in\extendedspace_{n+1} \Big\vert \kappa= {
\sum_{h\in\extendedspace_{n+1}} m_h} \biggr\}.\nonumber
\end{eqnarray}
Because $\zeta_1,\zeta_2,\ldots$ are i.i.d.,
conditioned on the sum $\kappa$,
the collection $\gprocess{\historycount{h}}{ h\in\extendedspace
_{n+1} }$
has a multinomial distribution, and so
%
%e47 #&#
\begin{eqnarray}
&&\mathbb{P}\biggl\{ \historycount{h} = m_h \colon h \in
\extendedspace_{n+1} \Big\vert \kappa= { \sum
_{h\in\extendedspace_{n+1}} m_h } \biggr\}
\nonumber
\\[-8pt]
\\[-8pt]
\nonumber
&&\quad = \frac{(\sum_{h\in\extendedspace_{n+1}} m_h) !} { \prod_{h\in
\extendedspace_{n+1}} (m_h !) } \prod_{h\in\extendedspace_{n+1}} \bigl[
\digammadist\bigl( h(n+1) ; r, c+nr \bigr)^{m_h} \bigr].
\end{eqnarray}
It follows that
%
%e48 #&#
\begin{eqnarray}
&& \mathbb{P}\{ \historycount{h} = m_h \colon h\in
\Histories_{n+1} \vert \historycount{\oh} = m_{\oh} \colon\oh\in
\Histories_n \}\nonumber
\\
&& \quad =\frac{ (c T\harmonicsum{r,c+nr})^{\sum_{h\in
\extendedspace
_{n+1}}     m_h} }{ (\sum_{h\in\extendedspace_{n+1}} m_h) ! } \exp( - c T\harmonicsum{r,c+nr} )
\nonumber
\\[-8pt]
\label{eqconddist}
\\[-8pt]
\nonumber
&& \qquad{}\times \frac{ \prod_{\oh\in\Histories_n} (m_{\oh} !) }{ \prod_{h\in
\existingspace_{n+1}} (m_h !) } \prod_{h\in\existingspace_{n+1}}
\bigl[ \bnbdist\bigl( h(n+1) ; r, \historysum{h} - h(n+1), c + nr
\bigr)^{m_h} \bigr]
\\
&&\qquad{}\times \frac{(\sum_{h\in\extendedspace_{n+1}} m_h) !} { \prod_{h\in
\extendedspace_{n+1}} (m_h !) } \prod_{h\in\extendedspace_{n+1}}
\bigl[ \digammadist\bigl( h(n+1) ; r, c+nr \bigr)^{m_h} \bigr].
\nonumber
\end{eqnarray}

\begin{pf*}{Proof of Theorem~\protect\ref{resultNBIBPPMF}}
The proof is by induction.
The p.m.f. $\mathbb{P}\{ \historycount{h} = m_h \colon h\in
\Histories_1 \}
$ is
given by equation (\ref{eqdistX1}), which agrees with equation (\ref{eqNBIBPPMF}) for
the case $n=1$.
The conditional p.m.f. $\mathbb{P}\{ \historycount{h} = m_h \colon
h\in
\Histories_{n+1} \vert \historycount{\oh} = m_{\oh} \colon\oh\in
\Histories_n \}$ is given by equation (\ref{eqconddist}).
By the inductive hypothesis, the p.m.f. $\mathbb{P}\{ \historycount
{\oh} =
m_{\oh} \colon\oh\in\Histories_n \}$ is given by equation (\ref{eqNBIBPPMF}).
Then by equation (\ref{eqJointCombStruct}), we have
%
%e49 #&#
%e50 #&#
%e51 #&#
\begin{eqnarray}
&&\mathbb{P}\{ \historycount{h} = m_h
\colon h\in\Histories_{n+1} \}
\nonumber\\
&&\quad= \frac{ (cT)^{ (\sum_{\oh\in\Histories_n}   m_{\oh}) }
(c T\harmonicsum{r,c+nr})^{ (\sum_{h\in\extendedspace
_{n+1}}     m_h) } }{
\prod_{h \in\existingspace_{n+1} } (m_h !)
 \prod_{h\in\extendedspace_{n+1}} (m_h !) } \exp \Biggl(   - c T\sum
_{j=1}^{n+1} \harmonicsum {r,c+(j-1)r} \Biggr)
\nonumber\\
\label{eqfirstmarginal}
&&\qquad{}\times \prod_{h\in\existingspace_{n+1}} \Biggl[ \betafn
\bigl(\historysum{h} - h(n+1), c+nr\bigr) \prod_{j=1}^n
  \frac{ (r)_{h(j)} }{ h(j) ! }
\\
&&\qquad {}\times \bnbdist\bigl( h(n+1) ; r, \historysum{h} - h(n+1), c +
nr \bigr) \Biggr] ^{m_h}\nonumber
\\
&&\qquad{}\times \prod_{h\in\extendedspace_{n+1}} \bigl[ \digammadist
\bigl( h(n+1); r, c+nr \bigr) \bigr]^{m_h}.\nonumber
\end{eqnarray}
In the first product term on the right-hand side of equation (\ref{eqfirstmarginal}), note that, for every $h\in\existingspace_{n+1}$,
\begin{eqnarray*}
&&\betafn\bigl(\historysum{h} - h(n+1), c+nr\bigr) \prod
_{j=1}^n   \frac{ (r)_{h(j)} }{ h(j) ! }
 \bnbdist\bigl( h(n+1) ; r, \historysum{h} - h(n+1), c +
nr \bigr)
\\
&&\quad = \betafn\bigl(\historysum{h}, c+(n+1)r\bigr) \prod
_{j=1}^{n+1} \frac{ (r)_{h(j)} }{ h(j) ! }.
\end{eqnarray*}
In the second product term, note that
\begin{eqnarray}
&& \prod_{h\in\extendedspace_{n+1}} \bigl[ \digammadist\bigl( h(n+1); r,
c+nr \bigr) \bigr]^{m_h}
\nonumber
\\
&&\quad= \prod_{h\in\extendedspace_{n+1}} \biggl[ \harmonicsum{r,c+nr}^{-1}
\frac{ (r)_{h(n+1)} }{ h(n+1) ! } \betafn\bigl( h(n+1), c+(n+1)r\bigr) \biggr]^{m_h}
\nonumber
\\
&&\quad= \harmonicsum{r,c+nr}^{- (\sum_{h\in\extendedspace_{n+1}}     m_h
) } \prod_{h\in\extendedspace_{n+1}}
\Biggl[ \betafn\bigl( h(n+1), c+(n+1)r\bigr) \prod_{j=1}^{n+1}
\frac{ (r)_{h(j)} }{ h(j) ! } \Biggr]^{m_h},
\nonumber
\end{eqnarray}
where for the last equality, we have used the fact that $h(j)=0$ for
every $j\le n$ and $h\in\extendedspace_{n+1}$.
Note that $\sum_{\oh\in\Histories_n} m_{\oh} + \sum_{ h\in
\extendedspace_{n+1} } m_h = \sum_{h\in\Histories_{n+1} } m_h$.
Then equation (\ref{eqfirstmarginal}) is equal to
%
%e52 #&#
\begin{eqnarray}
&&\frac{ (cT)^{ \sum_{h\in\Histories_{n+1}} m_h } }{
\prod_{h\in\Histories_{n+1}} (m_h !) } \exp \Biggl(   - c T\sum
_{j=1}^{n+1} \bigl[ \psi( c + j r ) - \psi\bigl( c +
(j-1)r \bigr) \bigr] \Biggr)
\nonumber
\\[-8pt]
\\[-8pt]
\nonumber
&&\qquad{}\times \prod_{h\in\Histories_{n+1}} \Biggl[ \betafn\bigl(
\historysum{h}, c+(n+1)r\bigr) \prod_{j=1}^{n+1}
\frac{ (r)_{h(j)} }{ h(j)! } \Biggr]^{m_h}.
\end{eqnarray}
Noting\vspace*{1pt} that $\sum_{j=1}^{n+1} [ \psi( c + j r ) - \psi( c + (j-1)r )
] = \psi( c +(n+1)r) - \psi(c)$, we obtain the expression in equation (\ref{eqNBIBPPMF}) for $n+1$, as desired.
\end{pf*}

By construction, equation (\ref{eqNBIBPPMF}) defines the finite-dimensional
marginal distributions of the stochastic process
$\gprocess{\historycount{h}}{h\in\countspace}$ with index set
$\countspace\defas\bigcup_{n\in{\mathbb{N}}} \Histories_n$.
The exchangeability result given by Theorem~\ref{resultNBIBPexch} then
follows from the exchangeability of the sequence $X_{[n]}$.

%s5 #&#
\section{Applications in Bayesian nonparametrics}
\label{secinference}

In Bayesian latent feature models,
we assume that there exists a latent set of \emph{features}
and that each data point possesses some (finite) subset of the features.
The features then determine the distribution of the observed data.
In a nonparametric setting,
exchangeable sequences of simple point processes can serve as models
for the latent sets of features.
Similarly, exchangeable sequences of point processes, like those that
can be constructed from
beta negative binomial processes,
can serve as models of latent \emph{multisets} of features.
In particular, atoms are features and their (integer-valued) masses
indicate multiplicity.
In this section, we develop posterior inference procedures for
exchangeable sequences of beta negative binomial processes.

%%%%%%%%%%%%%%%%%%%%%%
%s5.1 #&#
\subsection{Representations as random arrays/matrices}
\label{secarrays}
%%%%%%%%%%%%%%%%%%%%%%

A convenient way to represent the combinatorial structure of an
exchangeable sequence of point processes is via an array/matrix $W$ of
non-negative integers, where the rows correspond to point processes and
columns correspond to atoms appearing among the point processes.
Informally, given an enumeration of the set of all atoms appearing in
$X_{[n]}$, the entry $W_{i,j}$ associated with the $i$th row and
$j$th column is the multiplicity/mass of the atom labeled $j$ in the
$i$th point process~$X_i$.

More carefully,
fix $n\in{\mathbb{N}}$ and
let $\gprocess{\historycount{h}}{h\in\Histories_n}$ be the
combinatorial structure of a sequence $X_1,\ldots,X_n$ of conditionally
i.i.d. negative binomial processes, given a shared beta process base
measure with concentration parameter $c>0$ and non-atomic base measure
$\widetilde{B}_0$ of finite mass $T$.
Let $\kappa\defas\sum_{h\in\Histories_n} \historycount{h}$ be the
number of unique atoms among $X_{[n]}$.
Then $W$ is an $n \times\kappa$ array of non-negative integers such
that, for every $h\in\Histories_n$, there are exactly $\historycount
{h}$ columns of $W$ equal to~$h$, where $h$ is thought of as a
length-$n$ column vector.
Note that $W$ will have no columns when $\kappa= 0$.

All that remains is to order the columns of $W$.
Every total order on $\Histories_n$ induces a unique ordering of the
columns of $W$.
Titsias \cite{Titsias2007} defined a unique ordering in this way, analogous
to the
left-ordered form defined by Griffiths and Ghahramani \cite{GG06} for the IBP.
In particular, for $h,h' \in\Histories_n$, let $\preceq$ denote the
lexicographic order given by: $h \preceq h'$ if and only if $h=h'$ or
$h(\eta) < h'(\eta)$, where $\eta$ is the first coordinate where $h$
and $h'$ differ.
We say $W$ is \emph{left-ordered} when its columns are ordered
according to $\preceq$.
Because there is a bijection between combinatorial structures
$\gprocess
{\historycount{h}}{h\in\Histories_n}$ and their unique representations
by left-ordered arrays,
the probability mass function of $W$ is given by equation (\ref{eqNBIBPPMF}).

Other orderings have been introduced in the literature:
If we permute the columns of $W$ uniformly at random,
then $W$ is the analogue of
the \emph{uniform random labeling} scheme described by
Broderick, Pitman and Jordan
\cite{broderick2013feature} for the IBP.
Note that the number of distinct ways of ordering the $\kappa$ columns
is given by the multinomial coefficient
%
%e53 #&#
\begin{eqnarray}\label{eqarrayorder}
&& \frac{ \kappa! }{ \prod_{h\in\Histories_n} \historycount{h} ! },
\end{eqnarray}
where the denominator arises from the fact that there are
$\historycount
{h}$ indistinguishable columns for every history $h \in\Histories_n$.
The following result is then immediate:
%
%%%%%%%%%%%%%%%%%%%%%%%%%%%%%%%%%%%%%%%%%%%%%%%%%%%%%
%
%th6 #&#
\begin{thm} \label{eqarraypmf}
Let $W$ be a uniform random labeling of $\gprocess{\historycount
{h}}{h\in\Histories_n}$ described above,
let $w \in\mathbb{Z}_+^{n\times k}$ be an array of non-negative integers
with $n$ rows and $k\ge0$ non-zero columns, and for every $j\le k$,
let $s_j \defas\sum_{i=1}^n w_{i,j}$ be the sum of column $j$.
Then
%
%e54 #&#
\begin{eqnarray}
&&\mathbb{P}\{ W  =  w \} = \frac{ (c T)^k }{ k! } \exp \bigl(   - c T\bigl[ \psi(
c + n r ) - \psi( c ) \bigr] \bigr) \prod_{j=1}^k
\Biggl[ \frac{ \Gamma(s_j) \Gamma(c+nr) }{ \Gamma( s_j + c + nr ) } \prod_{i=1}^n
  \frac{ (r)_{w_{i,j}} }{ w_{i,j} ! } \Biggr].\qquad\hspace*{6pt}
\end{eqnarray}
\end{thm}
%
%%%%%%%%%%%%%%%%%%%%%%%%%%%%%%%%%%%%%%%%%%%%%%%%%%%%%

%f1
%f1 #&#
\begin{figure}[t]

\includegraphics{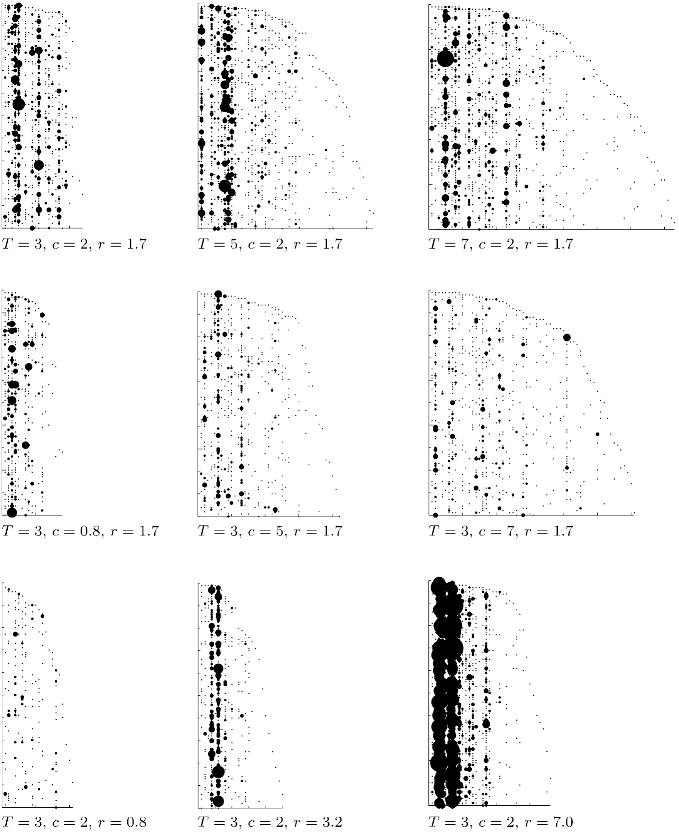}

\caption{Simulated $\mathbb{Z}_+$-valued arrays from the NB-IBP. Dots are
positive entries, the magnitudes of which determine the size of the
dot. The total mass parameter $T$ is varied along the top row;
the concentration parameter $c$ is varied along the middle row;
the negative binomial parameter $r$ is varied along the bottom row.
See the text for a summary of how these parameters affect the expected
number of features in total, features per row, and feature multiplicities.}
\label{figpriorsimulations}
\end{figure}

An array representation makes it easy to visualize some properties of
the model.
For example, in Figure~\ref{figpriorsimulations} we display several
simulations from the NB-IBP with varying values of the parameters
$T, c$, and $r$.
The columns are displayed in the order of first appearance, and are
otherwise ordered uniformly at random. (A similar ordering was used by
Griffiths and Ghahramani
\cite{GG06} to introduce the IBP.)
The relationship of the model to the values of $T$ and $c$ are
similar to the characteristics described
by
Ghahramani,
Griffiths and
Sollich
\cite{GGS2007} for the
IBP, with the parameter $r$ providing flexibility with respect to the
counts in the array.
In particular, the total number of features, $\kappa$, is Poisson
distributed with mean $c T[ \psi( c + n r ) - \psi( c ) ]$,
which increases with $T$, $c$, and $r$.
From the NB-IBP, we know that the expected number of features for the
first (and therefore, by exchangeability, every) row is $T$.
Because of the ordering we have chosen here, the rows are not
exchangeable, despite the sequence $X_{[n]}$ being exchangeable. (In
contrast, a uniform random labeling $W$ is row exchangeable and,
conditioned on $\kappa$, column exchangeable.)
Finally, note that the mean of the $\digammadist(r, c)$ distribution
exists for $c>1$ and is given by
%
%e55 #&#
\begin{eqnarray}
&& \frac{r} {(c-1) (\psi(r+c) - \psi(c))},
\end{eqnarray}
which increases with $r$ and decreases with $c$.
This is the expected multiplicity of each feature for the first row,
which again, by exchangeability, must hold for every row.
We may therefore summarize the effects of changing each of these
parameters (as we hold the others constant) as follows:
\begin{itemize}
\item Increasing the mass parameter $T$ increases both the
expected total number of features and the expected number of features
per row, while leaving the expected multiplicities of the features unchanged.

\item Increasing the concentration parameter $c$ increases the expected
total number of features and decreases the expected multiplicites of
the features, while leaving the expected number of features per row unchanged.

\item Increasing the parameter $r$ increases both the expected total
number of features and the expected multiplicities of the features,
while leaving the expected number of features per row unchanged.
\end{itemize}
These effects can be seen in the first, second, and third rows of
Figure~\ref{figpriorsimulations}, respectively.
We note that $r$ has a weak effect on the expected total number of
features (seen in the third row of Figure~\ref{figpriorsimulations}), and
$c$ has a weak effect on the expected multiplicities of the features
(seen in the second row of Figure~\ref{figpriorsimulations}).
The model may therefore be effectively tuned with $T$ and $c$
determining the size and density of the array, and $r$ determining the
multiplicities.
The most appropriate model depends on the application at hand, and in
Section~\ref{secconditionaldist} we discuss how these parameters may be
inferred from data.

%s5.2 #&#
\subsection{Examples}
\label{secexamples}

Latent feature models with associated multiplicities and unbounded
numbers of features
have found several applications in Bayesian nonparametric statistics,
and we now provide some examples.
In these applications, the features represent latent objects or factors
underlying a dataset
comprised of $n$ groups of measurements $y_1,\ldots,y_n$, where
each group $y_i$ is comprised of $D_i$ measurements $y_i =
(y_{i,1}, \ldots, y_{i,D_i})$.
In particular, $W_{i,j}$ denotes the number of instances of
object/factor $j$ in group $i$.

These nonparametric latent feature representations lend themselves
naturally to mixture models with an unbounded number of components.
For example, consider a variant of the models by
Sudderth \textit{et al.}
\cite{2005describing}
and
Titsias
\cite{Titsias2007} for a dataset of $n$ street camera images where
the latent features are interpreted as object classes that may appear
in the images, such as ``building'', ``car'', ``road'', etc.
The count $W_{i,j}$ models the relative number of times object class
$j$ appears in image $i$.
For every $i\le n$, image $y_i$ consists of $D_i$ local patches $y
_{j,1}, \ldots, y_{j,D_i}$ detected in the image, which are
(collections of) continuous variables representing, for example, color,
hue, location in the image, etc.
Let $\kappa$ be the number of columns of $W$, that is, the number of features.
The local patches in image $i$ are modeled as conditionally i.i.d. draws
from a mixture of $S_i = \sum_{j=1}^\kappa W_{i,j}$\vspace*{1pt} Gaussian distributions,
where $W_{i,j}$ of these\vspace*{1pt} components are associated with feature $j$.
For $k=1,2,\ldots,$ let
$\Theta_{i}^{(j,k)} \defas(m_{i}^{(j,k)}, \Sigma_{i}^{(j,k)})$
denote the mean and covariance of the Gaussian components associated
with feature $j$ for image $i$.
Let $z_{i,d} = (j,k)$ when $y_{i,d}$ is assigned to component $k \le W_{i,j}$
associated with feature $j \le\kappa$.
Conditioned on $\Theta\defas(\Theta_{i}^{(j,k)})_{i\le n, j \le
\kappa, k\le W_{i,j}}$ and the assignments $z \defas(z_{i,d})_{i\le
n, d\le D_i}$,
the distribution of the measurements admits a conditional\vspace*{-3pt} density
%
%e56 #&#
\begin{eqnarray}
\label{eqlikel1}
&& p(y \vert W, \Theta, z ) = \prod_{i=1}^n
\prod_{d=1}^{D_i} \mathcal{N} \bigl(
y_{i,d}; m_{i}^{z_{i,d}}, \Sigma_{i}^{z_{i,d}}
\bigr).
\end{eqnarray}
To share statistical strength across images, the parameters $\Theta
_{i}^{(j,k)}$ are given a hierarchical Bayesian prior:
%
%e57 #&#
%e58 #&#
\begin{eqnarray}
\Theta_{i}^{(j,k)} \vert \Theta^{(j)} &\stackrel{\mathrm{i.i.d.}}
{\sim} & \nu\bigl(\Theta^{(j)}\bigr)\qquad \mbox{for every $i$ and $k$,}
\\
\Theta^{(j)} &\stackrel{\mathrm{i.i.d.}} {\sim} & \nu_0\qquad \mbox{for
every $j$.}
\end{eqnarray}
A typical choice for
$\nu(\cdot)$
is the family of Gaussian--inverse-Wishart distributions with
feature-specific parameters $\Theta^{(j)}$ drawn i.i.d. from a
distribution $\nu_0$.
Finally, for every image $i\le n$, conditioned on $W$, the assignment
variables $z_{i,1}, \ldots, z_{i,D_i}$ for the local patches in image
$n$ are
assumed to form a multivariate P\'olya urn scheme, arising from
repeated draws from a Dirichlet-distributed probability vector over $\{
(j,k) :j \le\kappa,  k \le W_{i,j} \}$. The parameters for the
Dirichlet distributions are tied in a similar fashion to $\Theta$.
The interpretation here is that local patch $d$ in image $i$ is
assigned to one of the $S_i$ instances of the latent objects appearing
in the image.
The number of object instances to which a patch may be assigned is
specific to the image, but components across all images that correspond
to the same feature will be similar.

Latent feature representations are also a natural choice for factor
analysis models.
Canny \cite{canny2004gap} and  Zhou \textit{et~al.} \cite{ZHDC12} proposed models for text
documents in terms of latent features representing \emph{topics}.
More carefully, let $y_{i,v}$ be the number of occurrences of word
$v$ in document $i$.
Conditioned on $W$ and a collection of non-negative topic-word weights
$\Theta\defas(\theta_{j,v})_{j\le\kappa, v\le V}$, the word counts
are assumed to be conditionally i.i.d. and
%
%e59 #&#
\begin{eqnarray}
\label{eqlikel2}
&& y_{i,v} \vert W, \Theta\sim \poisson \Biggl( \sum
_{j=1}^\kappa W_{i,j}
\theta_{j,v} \Biggr).
\end{eqnarray}
In other words, the expected number of occurrences of word $v$ in
document $i$ is a linear sum of a small number of weighted factors. The
features here are interpreted as topics: words $v$ such that $\theta
_{j,v}$ is large are likely to appear many times.
There are a total of $\kappa$ topics that are shared across the documents.
The topic-word weights $\Theta$ are typically chosen to be i.i.d. Gamma
random variates, although there may be reason to prefer priors with
dependency enforcing further sparsity.
This general setup has been applied to other types of data including,
for example, recommendations \cite{gopalan2014bayesian},
where $y_{i,v}$ represents the rating a Netflix user $i$ assigns to a
film $v$.

%%%%%%%%%%%%%%%%%%%%%%
%s5.3 #&#
\subsection{Conditional distributions}
\label{secconditionaldist}
%%%%%%%%%%%%%%%%%%%%%%

Let $W$ be a uniform random labelling of a NB-IBP as described in
Section~\ref{secarrays}.
In the applications described above, computing the posterior
distribution of $W$ is the first step towards most other inferential goals.
Existing inference schemes use stick-breaking representations, that is,
they represent (a truncation of) the beta process underlying $W$.
This approach has some advantages, including that the entries of $W$
are then conditionally independent negative binomial random variables.
On the other hand, the random variables representing the truncated beta
process, as well as the truncation level itself, must be marginalized
away using auxiliary variable methods or other techniques \cite{BMPJpre,ZHDC12,TGG07,PZWGC2010}.
Here, we take advantage of the structure of the NB-IBP and do not
represent the beta process. The result is a set of Markov (proposal)
kernels analogous to those originally derived for the IBP \cite{GG06,GGS2007}.

The models described in Section~\ref{secexamples} associate every feature
with a latent parameter. Therefore,
conditioned on the number of columns $\kappa$, let $\Theta= (\theta_1,
\ldots, \theta_\kappa)$ be an i.i.d. sequence drawn from some non-atomic
distribution $\nu_\Theta$,
and assume that the data $y$ admits a conditional density $p(y \vert W, \Theta)$.
We will associate the $j$th column of $W$ with $\Theta_j$, and so the
pair $(W,\Theta)$ can be seen as an alternative representation for an
exchangeable sequence $X_{[n]}$ of beta negative binomial processes.
By Bayes' rule, the posterior distributions admits a conditional density
%
%e60 #&#
\begin{eqnarray}
&& p(W,\Theta| y) \propto p ( y \vert W, \Theta) \times p ( W, \Theta),
\end{eqnarray}
where $p(W,\Theta)$ is a density for the joint distribution of
$(W,\Theta)$.
We describe two Markov kernels that leave this distribution invariant.
Combined, these kernels give a Markov chain Monte Carlo (MCMC)
inference procedure for the desired posterior.

The first kernel resamples individual elements $W_{i,j}$, conditioned
on the remaining elements of the array (collectively denoted by
$W_{-(i,j)}$), the data $y$, and the parameters $\Theta$.
By Bayes' rule, and the independence of $\Theta$ and $W$ given $\kappa
$, we have
%
%e61 #&#
\begin{eqnarray}
&& \mathbb{P} \{ W_{i,j} = z \vert y,
W_{-(i,j)}, \Theta\}
\nonumber
\\[-8pt]
\label{eqposteriorWgeneral}
\\[-8pt]
\nonumber
&&\quad\propto p \bigl( y \vert \{ W_{i,j} = z \}, W_{-(i,j)},
\Theta\bigr) \times \mathbb{P}\{ W_{i,j} = z \vert W_{-(i,j)} \}.
\end{eqnarray}

Recall that the array $W$ is row-exchangeable,
and so, in the language of the NB-IBP, we may associate the $i$th row
with the final customer at the buffet.
The count $W_{i,j}$ is the number of servings the customer takes of
dish $j$, which has been served $S_{j}^{(-i)} \defas\sum_{i' \ne i}
W_{i',j}$ times previously.
When $S_{j}^{(-i)} > 0$, we have
%
%e62 #&#
\begin{eqnarray} \label{eqconditionalW}
&& W_{i,j} \vert W_{-(i,j)} %, \{S_{j}^{(-i)} > 0\}
\sim\bnbdist\bigl( r,
S_{j}^{(-i)}, c + (n-1)r \bigr).
\end{eqnarray}
Therefore, we can simulate from the unnormalized, unbounded discrete
distribution in equation~(\ref{eqposteriorWgeneral}) using equation (\ref{eqconditionalW}) as a  Metropolis--Hastings proposal, or we could use
inverse transform sampling where the normalization constant is
approximated by an importance sampling estimate.

Following Meeds \textit{et al.} \cite{meeds2007modeling},
the second kernel resamples the number, positions, and values of those
\emph{singleton} columns $j'$ such that $S_{j'}^{(-i)} = 0$.
Simultaneously, we propose a corresponding change to the sequence of
latent parameters $\Theta$, preserving the relative ordering with the
columns of $W$. This corresponding change to $\Theta$ cancels out the
effect of the $\kappa!$ term appearing in the p.m.f. of the array $W$.
Let $J_i$ be the number of singleton columns, that is, let
%
%e63 #&#
\begin{eqnarray}
&& J_i = \cardsharp\bigl\{ j \le\kappa\colon W_{i,j} > 0
\mbox{ and } S_{j}^{(-i)} = 0 \bigr\},
\end{eqnarray}
which we note may be equal to zero.
Because we are treating the customer associated with row $i$ as the
final customer at the buffet,
$J_i$ may be interpreted as the number of new dishes sampled by the
final customer, in which case,
we know that
%
%e64 #&#
\begin{eqnarray}
\label{eqsingletondishdist}
&& J_i \sim\poisson\bigl( cT\bigl[ \psi(c+nr) - \psi
\bigl(c+(n-1)r\bigr) \bigr] \bigr).
\end{eqnarray}
We therefore propose a new array $W^*$ by removing the $J_i$ singleton
columns from the array and insert $J_i^*$ new singleton columns at
positions drawn uniformly at random, where $J_i^*$ is sampled from the
(marginal) distribution of $J_i$ given in equation (\ref{eqsingletondishdist}).
Like those columns that were removed, each new column has exactly one
non-zero entry in the $i$th row:
We draw each non-zero entry independently and identically from a
$\digammadist(r, c+(n-1)r)$ distribution, which matches the
distribution of the number of servings the last customer takes of each
newly sampled dish.

Finally, we form a new sequence of latent parameters $\Theta^*$ by
removing those entries from $\Theta$ associated with the $J_i$ columns
that were removed from $W$ and inserting $J_i^*$ new entries, drawn
i.i.d. from $\nu_\Theta$, at the same locations corresponding to the
$J_i^*$ newly introduced columns.
Let $\kappa^* \defas\kappa- J_i + J_i^*$,
and note that there were
${\kappa^* \choose J_i^*}$ possible ways to insert the new columns.
Therefore, the proposal density is
%
%e65 #&#
\begin{eqnarray}
q\bigl( W^*, \Theta^* \vert W, \Theta\bigr) &=&
\pmatrix{\kappa^*
\cr
J_i^*}^{-1} \poisson\bigl(
J_i^* ; c T\bigl[ \psi(c+nr) - \psi\bigl(c+(n-1)r\bigr) \bigr]
\bigr)
\nonumber
\\[-8pt]
\label{eqproposaldist}
\\[-8pt]
\nonumber
&&{}\times \prod_{j\le\kappa^*} \digammadist\bigl(
W_{i,j}^* ; r, c+(n-1)r \bigr) \prod_{\theta\in\Theta^* \setminus\Theta}
\nu_\Theta( \theta).
\end{eqnarray}
With manipulations similar to those in the proof of Theorem~\ref{resultNBIBPPMF}, it is straightforward to show that
a Metropolis--Hastings kernel accepts a proposal $(W^*, \Theta^*)$
with probability $\min\{ 1, \alpha^* \}$, where
%
%e66 #&#
\begin{eqnarray}
&&\alpha^* = \frac{ p(y \vert W^*, \Theta^*) }{ p(y \vert W, \Theta) }.
\end{eqnarray}
Combined with appropriate Metropolis--Hastings moves that shuffle the
columns of $W$ and resample the latent parameters $\Theta$,
we obtain a Markov chain whose stationary distribution is the
conditional distribution of $W$ and $\Theta$ given the data $y$.

Another benefit of the characterization of the distribution of $W$ in
(\ref{eqarraypmf}) is that numerically integrating over the real-valued
concentration, mass, and negative binomial parameters $c$, $T
$, and $r$, respectively, are straightforward with techniques such as
slice sampling \cite{Neal2003}.
In the particular case when $T$ is given a gamma prior
distribution, say $T\sim\gammadist(\alpha,\beta)$ for some
positive parameters $\alpha$ and $\beta$, the conditional distribution
again falls into the class of gamma distributions.
In particular,
the conditional density is
%
%e67 #&#
%e68 #&#
\begin{eqnarray}
p(T\vert W, \kappa) &\propto & T^{\alpha+ \kappa-1} \exp \bigl(   - c T\bigl[ \psi( c
+ n r ) - \psi( c ) \bigr] - \beta T \bigr)
\\
&\propto&  \gammadist \bigl( T; \alpha+ \kappa, \beta+ c T \bigl[ \psi( c + n r )
- \psi( c ) \bigr] \bigr).
\end{eqnarray}

%%%%%%%%%%%%%%%%%%%%%%
\section*{Acknowledgements}
%%%%%%%%%%%%%%%%%%%%%%

We thank Mingyuan Zhou for helpful feedback and for pointing out the
relation of our work to that of Sibuya \cite{Sibuya1979}. We also thank Yarin
Gal and anonymous reviewers for feedback on drafts.
This research was carried out while C. Heaukulani was supported by the
Stephen Thomas studentship at Queens' College, Cambridge, with funding
also from the Cambridge Trusts, and while D.M.~Roy was a research
fellow of Emmanuel College, Cambridge, with funding also from a Newton
International Fellowship through the Royal Society.

% imsref loaded by daiva.urboniene, 2015-08-14 15:49:53
% imsref loaded by giedrek, 2015-09-29 16:08:27
% imsref loaded by giedrek, 2015-09-29 16:09:16

\printhistory
\end{document}